# On optimal cover and its possible shape for fractals embedded into 2D Euclidian space

Dmitry Zhabin[1]

Integrated Risk Management Department

VTB Bank (PJSC), Moscow, 123112, Presnenskaya emb., 12

In this article a definition of optimal cover for fractal structures is proposed. Minkowsky dimension is rewritten as functional equation on areas of covers that constructed for different scales. Given the definition, the functional equation is resolved and possible shapes of optimal coverage are defined in correspondence with fractal dimension values.

## 1. Introduction.

Fractal dimension is supposed to reflect a correspondence between geometrical shape of a figure and a measure of the space in which the figure is embedded [1, 2]. Generally, topological dimension is an integer-valued number that has natural association with length, area and volumes depending on the dimension value. That is evidently true for Euclidian spaces and smooth geometrical shapes. So, fractal dimensions, which commonly can be non-integer number, can be treated as a measure of complexity of the figure's shape [3].

Any set embedded into Euclidian space can be covered by number of ordinary sets, which shape, form and area are well defined. Such collection of ordinary sets is denoted as a set cover. Sum of areas of those ordinary sets is an area of the cover. For defined size of those ordinary sets one can found such a way of cover making, which has minimal total area. Any change in the size of the ordinary sets triggers change in number of ordinary sets, which required for minimal cover. The dependence of number of ordinary sets for minimal cover on the size of those sets includes information about complexity of the embedded set shape. Commonly, through investigation of the dependence one can find the value of fractal dimension, although it may be a consuming calculation.

---

[1] E-mail address: dzhabin@vtb.ru

In this article we consider fractals embedded into 2D Euclidian space. The Minkowsky dimension is rewritten in terms of functional equation that expresses a correspondence between areas of covers constructed for different scales. It is shown that the functional equation has at least one solution. The solution found is based on definition of optimal cover and shows a link between shape of the optimal cover and fractal dimension, which is a measure of fractal irregularity. Finally, shapes of optimal cover are calculated for different values of fractal dimension.

## 2. Fractal dimension and functional equation.

As a definition of fractal dimension, the Minkowsky dimension is widely known [4], when an ordinary sets, which are used as a set cover, has simple geometrical shape like box, circle or similar:

$$\dim F = -\lim_{\delta \to 0} \ln(N(\delta))/\ln(\delta). \qquad (1)$$

Here $\delta$ - is a size of ordinary set; $N(\delta)$ – number of ordinary sets with size $\delta$ required to cover the set. Minkowsky dimension is not the only way to calculate fractal dimension, for more information about existing methods of fractal dimension estimation please see [5].

While expression (1) shows the correspondence between number of ordinary sets, its size and fractal dimention, one may ask a question what kind of rule is imposed on areas of sets cover for fractals? To get close to the answer on this question, let us rewrite the (1) in terms of areas of set cover rather than its numbers.

Let us assume further that the limit (1) exists and equal to $D_H$, then following identity is true:

$$N(\delta) = C \cdot \delta^{-D_H}, \qquad (2)$$

where $C$ – some positive constant.

For the sake of simplicity, let us assume further on that the fractal set under consideration is embedded into Euclidian space. By ordinary set let us assume a set of simple geometrical shape, which area is proportional to size $\delta$ of the set in integer power that is equal to dimension of the Euclidian space, like box, square, circle of sphere. So, the area of the set cover follows rule:

$$S(\delta) \sim N(\delta)\delta^{D_E}, \qquad (3)$$

where $D_E$ – is topological dimension of the Euclidian space; $\delta$ - is length (or radius) of the ordinary set.

If length of the ordinary set is changed in $\alpha$ times so that keeping $\alpha \cdot \delta \sim \delta$, i.e. $\alpha \cdot \delta$ has the same order as $\delta$, number of sets with length $\alpha \cdot \delta$ required to cover the set can be estimated as:

$$N(\alpha\delta) = C \cdot (\alpha\delta)^{-D_H} = \alpha^{-D_H} N(\delta). \qquad (4)$$

In the same time the area of the cover changes as:

$$S(\alpha\delta) \sim N(\alpha\delta)(\alpha\delta)^{D_E}. \qquad (5)$$

Combining (3) and (4) one can get the functional equation on the areas of the covers of different scales:

$$S(\alpha\delta) = \alpha^{D_E - D_H} S(\delta). \qquad (6)$$

Expression (6) shows the rule we sought at the beginning of this section. For any set embedded into Euclidian space, and for which limit (1) exists, areas of set cover follows rule (6) when a size of the ordinary set $\delta$ changes in $\alpha$ times.

From the other hand, expression (6) is a functional equation and a solution of this equation may give some extra information regarding possible shape of set cover that corresponds to complex nature of fractals mostly.

## 3. Optimal cover and its shape for 2D Euclidian case.

Usually one defines shape of ordinary set used for cover and then investigates the rule by which the number of ordinary sets depends on size of ordinary set. Here we will try to found the shape of the set used for cover when the rule by which the number of those sets depends on the set size is defined. For the clarity, let us in this section consider 2D Euclidian space.

**Definition.**

*Let's denote as an optimal set cover the one, for which number of covering sets required to get minimal coverage, changes by linear rule as the size of the covering set is scaled.*

That means that when optimal cover covers the fractal set under consideration, then if the size of the covering sets is scaled, so the number of covering sets changes linearly. So, expression for area of set cover can be defined:

$$S(\delta) = \sum_{i=1}^{N'(\delta)} \varphi(\delta) = N'(\delta)\varphi(\delta), \qquad (7)$$

where $\varphi(\delta)$ – is area of single optimal set; $N'(\delta)$ – number of sets in optimal set cover. So, if for 2D Euclidian case the size of the optimal set is scaled then area of the optimal cover changes as:

$$S(\alpha\delta) = \sum_{i=1}^{N'(\alpha\delta)} \varphi(\alpha\delta) = N'(\alpha\delta)\varphi(\alpha\delta) = \frac{1}{\alpha} N'(\delta)\varphi(\alpha\delta). \qquad (8)$$

From the other hand, areas of set covers are to obey expression (6), and then one can easily get functional equation on areas for single optimal set:

$$\varphi(\alpha\delta) = \alpha^{D_E - (D_H - 1)} \varphi(\delta) . \qquad (9)$$

Thus, if the optimal set cover exists it is to satisfy equation (9). Solution of (9) can be a unique source of information how nontrivial fractal geometry is reflected in the shape of set cover.

It can be shown, see [6], that following function satisfies (9):

$$\varphi(\delta) = C_1 \cdot \delta^{D_E - (D_H - 1)}, \qquad (10)$$

where $C_1$ – is some positive constant.

Having (10) one can find shape of the optimal set cover. Indeed, for 2D Euclidian space we can assume that single optimal set cover is defined on interval $[0,\delta]$, so the shape of the border of optimal cover can be found on class of monotonic functions:

$$\varphi(\delta) = C_1 \cdot \delta^{D_E - (D_H - 1)} = C_2 \int_0^{\delta} f(x) \, dx, \qquad (11)$$

where $f(x) = \begin{cases} 0, & \text{if } x \notin [0,\delta] \\ x^{D_E - D_H}, & \text{if } x \in [0,\delta] \end{cases}$

Expression (11) allows a geometrical interpretation, see Picture 1, where shape of single optimal set cover is shown for different cases of fractal dimensions.

It is easy to see that for $D_H = 1$ case, which is typical for smooth curves, optimal cover is an ensemble of points bounded by triangle. As the fractal dimension growths the border of optimal cover behaves more irregularly and the close fractal dimension to space dimension the close shape of optimal cover to an rectangle.

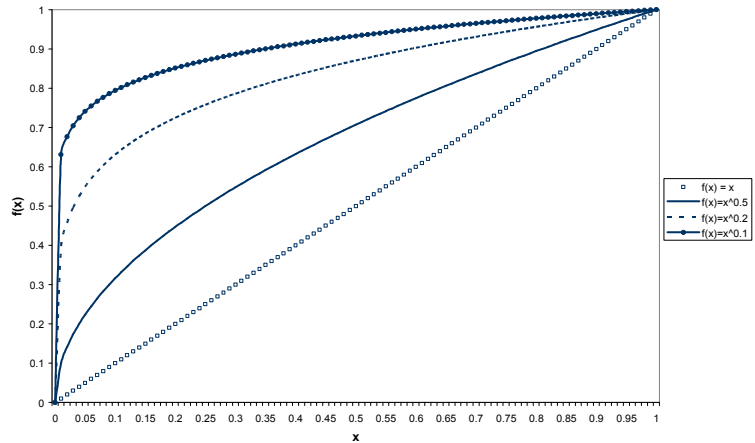

**Picture 1.** *Dependence of shape of optimal set cover from difference between space and fractal dimensions.*

. So, it is shown that at least one solution of equation (6) exists, this solution is based on definition of optimal set cover and a shape of this optimal cover differs for different values of fractal dimensions.

## Conclusions

Fractal dimension is usually treated as a proper measure of the system complexity. Commonly, this complexity may be found in dependence of number of ordinary sets on the scale of those ordinary sets. Equivalently, this complexity imposes the functional equation (6) on areas of set covers.

Among all possible ways to cover fractals with set cover there are certain ways that allows solving the equation (6). This certain way of cover construction we denote as optimal cover and, in a way, it "transfers" irregular dependence of number of ordinary sets that used to cover a fractal to a shape of ordinary set. So, at least one solution (10) of functional equation (6) was found.

For 2D case the shape of optimal cover tends from triangle to rectangle as fractal dimension goes from 1 to 2. Such a behavior may have very natural interpretation for 2D time-series if applied to positive and negative autocorrelated series.

## References.